\documentclass[12pt,reqno]{article}
mm \textheight=8.9in

\usepackage{amssymb}
\usepackage{amsmath}
\usepackage{amsthm}

\usepackage{color}

\newcommand{\br}[3]{{$#1$}$\lower4pt\hbox{$\tp\atop\raise4pt \hbox{$\scriptscriptstyle{#2}$}$} ${$#3$}}
\newcommand{\tw}[3]{{$#1$}${\,\scriptscriptstyle {#2}}\atop\raise9pt\hbox{$\scriptstyle\tp$} ${$#3$}}
\newcommand{\ttps}[2]{{#1}\raise5pt\hbox{$\lower12pt\hbox{$\scriptstyle\tp$}\atop \lower0pt\hbox{$\tilde\;$}$}\raise4.5pt\hbox{${\scriptstyle{#2}}$}}
\newcommand{\st}[1]{\mbox{${\,\scriptscriptstyle {#1}}\atop\raise5.5pt\hbox{$*$}$}}
\newcommand{\btr}{\raise1.2pt\hbox{$\scriptstyle\blacktriangleright$}\hspace{2pt}}

\newcommand{\ttp}{{\lower12pt\hbox{$\tp$}\atop \hbox{$\tilde\;$}}}

\newcommand{\id}{\mathrm{id}}

\newcommand{\A}{\mathcal{A}}
\newcommand{\B}{\mathcal{B}}

\newcommand{\Der}{\mathrm{Der}}
\newcommand{\Lie}{\mathrm{Lie}\>}

\newcommand{\Ha}{\mathcal{H}}
\newcommand{\Ru}{\mathcal{R}}

\newcommand{\E}{\mathcal{E}}

\newcommand{\C}{\mathbb{C}}

\newcommand{\tp}{\otimes}
\newcommand{\vt}{\vartheta}

\newcommand{\zt}{\zeta}

\newcommand{\Ur}{\mathrm{U}}

\newcommand{\F}{\mathcal{F}}
\newcommand{\ve}{\varepsilon}
\newcommand{\gm}{\gamma}

\newcommand{\prt}{\partial}

\newcommand{\tr}{\triangleright}

\newcommand{\End}{\mathrm{End}}

\newcommand{\btl}{\mbox{\raise1.1pt\hbox{$\scriptstyle\labcktriangleleft$}}}
\newcommand{\ad}{\mathrm{ad}}

\newcommand{\g}{\mathfrak{g}}

\newcommand{\nn}{\nonumber}

\newcommand{\T}{\mathrm{T}}

\newcommand{\al}{\alpha}

\newcommand{\be}{\begin{eqnarray}}
\newcommand{\ee}{\end{eqnarray}}

\newtheorem{thm}{Theorem}[section]
\newtheorem{propn}[thm]{Proposition}

\newtheorem{corollary}[thm]{Corollary}

\theoremstyle{definition}
\newtheorem{remark}[thm]{Remark}
\newtheorem{definition}[thm]{Definition}
\newtheorem{example}[thm]{Example}

\newcount\prg

\newcommand{\parag}{\advance\prg by1 {\noindent\bf\thesection.\the\prg\hspace{6pt}}}

\begin{document}
\title{Twisting Adjoint Module Algebras}
\author{Petr Kulish
 \\
\small St.-Petersburg Department of Steklov Mathematical Institute,\\
\small Fontanka 27, 191023 St.-Petersburg, Russia
\and Andrey Mudrov \\
\small Department of Mathematics, University of Leicester, \\
\small University Road,
LE1 7RH Leicester, UK}
\date{}
\maketitle

\begin{abstract}
Transformation of operator algebras under Hopf algebra twist is studied.
It is shown that that adjoint module algebras are stable under the twist.
Applications to vector fields on non-commutative space-time are considered.
\end{abstract}



\section{Introduction}
The interest in the quantum field theory on non-commutative space-time
arose long before the invention of quantum groups. It was reinvigorated after
the Heisenberg relations for the coordinates appeared as a special limit in
the open string theory of Seiberg-Witten, \cite{SW}. Later those
relations were  shown to be invariant under an action of the quantum Poincar\'e  algebra,
which resulted from the classical Poincar\'e algebra through a quantum group twist \cite{CKNT}.

Twist transformation of the Poincar\'e algebra leads
to a deformation of algebraic structures on the space-time, such as the function algebra, differential operators {\em etc}. In this paper we study the behavior of differential operators under this transformation.
In general mathematical terms, we consider a Hopf algebra $\Ha$, an adjoint module algebra $\A$, and how
the latter changes under a twist of $\Ha$. Recall that adjoint module algebra is
an arbitrary algebra admitting a homomorphism $\Ha\to \A$, which defines an (adjoint) action of $\Ha$ on $\A$
compatible with the multiplication  in $\A$. Typical examples come from representations of $\Ha$: one can take  for $\A$
the image of $\Ha$, the entire  algebra of endomorphisms, or any algebra in between.
We prove that twist transformation of an adjoint module algebra $\A$ is isomorphic to $\A$.

Our interest in the adjoint action is motivated by geometrical applications of quantum groups.
We regard  $\Ha$ as symmetry of the quantum space $\A$, whose geometry also involves operators
acting on $\A$ (like differential operators in the differential setting). In this picture, the algebra $\Ha$ is a sort of thing that is external to the geometry of the quantum space. Indeed, under the twist  of $\Ha$ the entire geometry transforms accordingly, see e.g. \cite{KM,K}, while the  multiplication in $\Ha$ remains untouched.

On the other hand, the Hopf algebra $\Ha$ "participates" in geometry through  the representation $\rho\colon \Ha\to \End(\A)$.
For example, the d'Alembertian $\square =\prt_\mu\prt^\mu$ on functions on the Minkowski space belongs to the (image of) Poincar\'e algebra. Thus, the multiplication between the partial derivatives in $\prt_\mu\prt^\mu$
must be deformed in $\End(\A)$,  but at the same time it must remain the same in $\Ha$. This apparent controversy is a cause of certain
confusion in the literature on the non-commutative gauge field theory (c.f. \cite{CT,ADMSW}). It can be resolved in the following way.
Instead of "external" object $\Ha$ one should deal with the image $\rho(\Ha)$. Under the twist, the algebra  $\rho(\Ha)$ is undergone
a deformation, however isomorphic to the original algebra $\Ha$. In the "new coordinates" both comultiplication and multiplication are deformed in a compatible way, and the algebra structure in $\Ha$  is incorporated uniformly into the whole picture.

This point of view helps to avoid the confusion between the deformed and non-deformed products,
as it allows to get rid of the  "external" object $\Ha$ and work with the purely geometrical object $\rho(\Ha)$. In particular, the noncommutative star product of  annihilation operators is often interpreted in the physics literature as a change of statistics. In the quantum theory $\A$ can be identified with the algebra of differential operators on the space-time (or the algebra of observables) while $\Ha$ is the symmetry
(Poincar\'e) algebra. Then the above mentioned isomorphism means that one cannot speak about
the "change of statistics" under a twist of comultiplication of $\Ha$ and the
corresponding multiplication in $\A$ (see Example \ref{example}).

Here is the setup of the paper. Some facts about Hopf algebras and twists that are used in what follows are collected  in the next part of Introduction. The second section is devoted to adjoint module algebras and their twist transformation. There, we also consider a special case of smash product of a Hopf algebra and its general
module algebra. Further we apply the general considerations of this section to universal enveloping Lie algebras
and their generalizations associated with unitary quantum permutations. This is illustrated on the example of quantum vector fields in the last section.
 
\subsection{Triangular Hopf algebras}
For the reader's convenience we collect in this section the basic information on Hopf algebras,
which will be used further on.
Throughout the paper, $\Ha$ is a {\em triangular} Hopf algebra over $\C$,
\cite{Dr2}.
That is, $\Ha$ is a complex associative algebra equipped with  comultiplication
$\Delta\colon \Ha\to \Ha\tp \Ha$, counit $\ve\colon \Ha\to \C$,
antipode $\gm\colon\Ha\to \Ha$, and a universal R-matrix $\Ru\in \Ha\tp \Ha$ which
obeys the identity $\Ru\Ru_{21}=1\tp 1$ (the subscripts indicate that the second copy of the R-matrix has the tensor legs flipped).

The comultiplication is an algebra homomorphism satisfying the
coassociativity condition $(\Delta\tp \id)\circ \Delta=(\id\tp \Delta)\circ \Delta$.
We shall use the standard symbolic Sweedler notation
$\Delta(h)=h^{(1)}\tp h^{(2)}$ for the coproduct of $h\in \Ha$. This notation
assumes the suppressed summation over decomposable tensors.

The counit $\ve$ is an algebra homomorphism, and
$(\ve\tp \id)\circ \Delta=\id=(\id\tp \ve)\circ \Delta $,
under the identification $\C\tp \Ha\simeq \Ha\simeq \Ha\tp \C$.

The antipode is an anti-algebra  and anti-coalgebra map. The latter means
that $(\gm \tp \gm)\circ \Delta=\Delta^{op} \circ  \gm$, where
$\Delta^{op}(h):=h^{(2)}\tp h^{(1)}$ is the opposite comultiplication.
The antipode satisfies the equalities
$\gm(x^{(1)})x^{(2)}=\ve(x)=x^{(1)}\gm(x^{(2)})$.

The universal R-matrix $\Ru$ is an invertible element from $\Ha\tp \Ha$ (may be a completed tensor product)
satisfying the identities \cite{Dr1}
$$
(\Delta\tp \id)(\Ru)=\Ru_{13}\Ru_{23},
\quad
(\id\tp \Delta)(\Ru)=\Ru_{13}\Ru_{12},
$$
in $\Ha^{\tp 3}$. Here the subscripts indicate the embeddings of $\Ru$
in $\Ha^{\tp 3}$ in the standard way.
Also,
$$
\Ru\Delta(h)=\Delta^{op}(h)\Ru,
$$
for all $h\in \Ha$.
As $\Ha$ is assumed to be triangular, $\Ru_{21}=\Ru^{-1}$.

Note that $(\Ha,\Delta^{op},\ve,\gm^{-1},\Ru^{-1})$ is again a triangular Hopf
algebra. This makes sense, as the antipode $\gm$ is invertible  in triangular Hopf algebras.
We denote it by $\Ha^{op}$. One can also define
the triangular Hopf algebra $\Ha_{op}$ with the opposite multiplication,
the same comultiplication and counit, and the inverse antipode and R-matrix.

An example of triangular Hopf algebra is a  universal enveloping algebra $\Ur(\g)$
of a Lie algebra $\g$ with
$$
\Delta(\xi)=\xi\tp 1+1\tp \xi,
\quad \ve(\xi)=0,
\quad
\gm(\xi)=-\xi
$$
on the elements $\xi\in \g$. The universal R-matrix is just $1\tp 1$.

An associative algebra $\A$ is called an $\Ha$-module algebra if
it is an $\Ha$-module and the multiplication of $\A$ is
compatible with the comultiplication of $\Ha$:
$$
h\tr(ab)=(h^{(1)}\tr a)(h^{(2)}\tr b),
\quad h\in \Ha, \quad a,b\in \A.
$$
In geometrical applications, $\Ha$ acts on a (non-commutative) space, $M$,
and   $\A$ is the algebra of "smooth" functions on $M$.
Observe that for any $\Ha$-module $(V,\rho)$ the vector space
$\End(V)$ of all linear operators on $V$ is an
$\Ha$-module algebra with the action
$h\tr A:=\rho(h^{(1)})A\rho(\gm(h^{(2)}))$.

Now we recall the basics of the Hopf algebra twist.
Suppose there is an invertible element $\F\in \Ha\tp \Ha$
satisfying the identities \cite{Dr2}
$$
(\Delta\tp \id)(\F)\F_{12}=(\id \tp \Delta)(\F)\F_{23},
$$
$$
(\ve\tp \id)(\F)=1\tp 1 =(\id\tp \ve)(\F).
$$
Then $(\Ha,\tilde\Delta,\ve,\tilde \gm,\tilde\Ru)$ is again a triangular
Hopf algebra with
$$
\tilde \Delta(h)\colon =\F^{-1}\Delta(h)\F
\quad
\tilde\gm(g)\colon =\vt^{-1}\gm(h)\vt,
\quad
\tilde \Ru\colon =\F^{-1}_{21}\Ru\F
$$
for all $h\in \Ha$.
The invertible element  $\vt$ is defined by the formula (\ref{modului}) below.
To distinguish this new Hopf algebra from $\Ha$, we denote it by $\tilde \Ha$
and call it twist of $\Ha$ via the cocycle $\F$.
Note that $\tilde \Ha$ has the same multiplication and counit as $\Ha$.

A twist of $\Ha$ transforms any $\Ha$-module algebra $\A$ to
an $\tilde \Ha$-module algebra $\A_\wr$.
As an $\Ha$-module, $\A_\wr$ coincides with $\A$ but
has a different multiplication,
$a*b:=(\F_1\tr a)(\F_2\tr b)$. This multiplication is compatible
with the twisted comultiplication in $\tilde \Ha$.
Thus a twist transforms Hopf algebras and their module algebras.
To distinguish transformation of module algebras from Hopf algebras, we
call $\A_\wr$  {\em cotwist} of $\A$.
\section{Twisting adjoint module algebras}
In the present subsection we consider module algebras of a special type.
Suppose $\F\in \Ha\tp \Ha$ is a twisting cocycle.
Define the elements $\vt,\zt\in \Ha$
by setting
\be
\vt=\gm(\F_1)\F_2,
\quad
\zt= \F_2^{-1} \gm^{-1}(\F_1^{-1}).
\label{modului}
\ee
It is known, \cite{Dr2}, that
\be
\vt^{-1}=\gm(\zt),
\quad
\zt^{-1}=\gm^{-1} (\vt).
\label{eq001}
\ee
As usual, the subscripts are used to label  tensor factors in the symbolic notation $\F=\F_1\tp \F_2$.
The elements $\vt,\zt$ participate in the antipode $\tilde \gm$ of $\tilde \Ha$, which is related to the
 old antipode by the formulas
\be
\tilde\gm(h)=\vt^{-1} \gm(h)\vt=\gm(\zt^{-1} h\zt).
\ee
Also, we have the identities
\be
\F^{(1)}_1\tp \gm(\F^{(2)}_1)\F_2=\F^{-1}_1\tp \gm(\F^{-1}_2)\vt,
\quad
\gm(\F_1)\F^{(1)}_2\tp \F^{(2)}_2=\vt \F^{-1}_1\tp \F^{-1}_2.
\label{eq003}
\ee
These formulas can be easily derived from the definition $\F$.

Now suppose
an associative algebra $\E$ admits a  homomorphism $\rho \colon \Ha\to \E$.
The adjoint action  $\ad_\rho\colon \Ha \to \End(\E)$ is defined by
$$
\ad_\rho(h) a=\rho(h^{(1)})a\rho(\gm(h^{(2)})), \quad h\in \Ha, \quad a \in \E.
$$
When $\rho$ is clear from the context, we use for the adjoint action the dot notation $\ad_\rho(h) a=h.a$.
The adjoint action makes $\E$ an $\Ha$-module algebra, which is further called adjoint module algebra.
We also extend this shorthand notation  to the action on maps between
any $\Ha$-modules, say  $(S,\rho_S)$ and $(T,\rho_T)$:
$$
h. f=\rho_T\bigl(h^{(1)})f \rho_S\bigl(\gm(h^{(2)})\bigr),
$$
for $f\colon S\to T$.
This will help us to distinguish the action on maps from the actions on $S$ and $T$.

We are going to show that $\E$ is stable under twist.
To simplify the formulas, we shall drop the symbol  of the homomorphism $\rho$, as though $\Ha$ were a subalgebra in $\E$.
\begin{thm}
\label{twistadjoint}
Let $\F\in \Ha\tp \Ha$ be a twisted cocycle and suppose $\E$ is an adjoint $\Ha$-module algebra.
Then the cotwist $\E_\wr$  is isomorphic to $\E$ as an associative algebra, with the isomorphism $\E_\wr\to \E$ given by the assignment
$\varphi \colon a\mapsto \F^{-1}_1a\gm(\F^{-1}_2)\vt=\bigl(\ad(\F_1)a\bigr) \F_2$.
\end{thm}
\begin{proof}
Denote the inverse mapping $\E\to \E_\wr$ by
$\varphi^{-1}\colon a\mapsto \F_1a\gm(\F_2\zt)$ and check that $\varphi^{-1}$ is an algebra homomorphism.
This immediately follows from the identity
\be
\F^{(1)}_1\F_{1'}\tp \gm(\F^{(2)}_1\F_{2'}\zt)\F^{(1)}_2\F_{1''}\tp \gm(\F^{(2)}_2\F_{2''}\zt)
=
\F_1\tp 1\tp \gm(\F_2\zt),
\label{eq004}
\ee
so let (\ref{eq004}) be checked first.
Applying the cocycle equation to $(\Delta\tp \Delta)(\F)\F_{12}$ we transform the left-hand side
to
\be
\F_1\tp \gm(\F^{(1)}_2\F_{1'}\zt)\F^{(2)}_2\F^{(1)}_{2'}\F_{1''}\tp \gm(\F^{(3)}_2\F^{(2)}_{2'}\F_{2''}\zt),
\nn
\ee
from which we easily get the expression
\be
\F_1\tp \gm(\zt)\gm(\F_{1'})\F^{(1)}_{2'}\F_{1''}\tp \gm(\F_2\F^{(2)}_{2'}\F_{2''}\zt).
\nn
\ee
Applying the right formula from (\ref{eq003}) and the left formula from (\ref{eq001}),
we obtain the right-hand side of (\ref{eq004}).
\end{proof}

\begin{example}
\label{example}
The universal enveloping algebra $\Ur({\cal P})$ of the Poincar\'e Lie algebra ${\cal P}$
is realized in the algebra  $\A$ of differential operators on
the Minkowski space, which is generated by the mutually commuting
coordinates $x^{\mu}$ and constant vector fields $i \partial/\partial x_\mu$ representing
the momenta $P_{\nu}\in {\cal P}$. Simple Abelian twist
$$
{\cal F}_1 \otimes {\cal F}_2:= {\cal F}= \exp \Bigl(-{\frac{i}{2}}
\theta^{\mu \nu} P_\mu \otimes P_\nu \Bigr), \quad \theta^{\mu \nu}=-\theta^{\nu\mu},
$$
changes the coproduct of $\Ur({\cal P})$ to
$\tilde\Delta(h) = {\cal F}^{-1} \Delta(h) {\cal F} $, for $h\in \Ur({\cal P})$. The algebra $\A$ is
an adjoint $\Ur({\cal P})$-module algebra, and after the twist it has
the $*$-product of $\A_\wr$:
$$
x^{\alpha} * x^{\beta} - x^{\beta} * x^{\alpha} = [x^\alpha , x^\beta ]_* = i
\theta^{\alpha \beta}.
$$
The homomorphism $\varphi\colon \A_\wr\to \A$  from Theorem \ref{twistadjoint}
acts by
$$
x\mapsto y^{\alpha}\mapsto ({\cal F}_1 x^\alpha) {\cal F}_2 = x^\alpha +
{\frac12} \theta^{\alpha \beta} P_\beta \in \A.
$$
One can check that the elements $y^\al$ obey the relations $[y^\alpha , y^\beta ] = i
\theta^{\alpha \beta} $.

Similar construction can be done for realization of $\Ur({\cal P})$ in
the algebra, call it $\A$ again, generated by  the creation and annihilation operators  $a(p)^{+}, a(k)$
of the relativistic scalar field. The action of $P_\mu$ on $a(p)$ is
$P_\mu \tr a(p) = [P_\mu, a(p) ] = - p_\mu a(p)$ giving the $*$-product
$a(p) * a(k) = a(k) * a(p) \exp ( - i \theta^{\mu \nu} p_\mu  k_\nu )$.
The isomorphism  $\varphi\colon \A_\wr\to \A$ results in
$$
a(p) \mapsto b(p) := ({\cal F}_1 a(p)) {\cal F}_2 = a(p) \exp ({\frac{i}{2}} \theta^{\mu \nu} p_\mu \otimes P_\nu ) \in
\A,
$$
$b(p) b(k) = b(k) b(p) \exp ( - i \theta^{\mu \nu} p_\mu  k_\nu )$. Therefore, $\A_\wr$
is realized in the smash product of the original algebra and $\Ur({\cal P})$ by a change of variables.

\end{example}
The isomorphism $\varphi$ can be defined in the obvious way for any adjoint module, not only for $\Ha$. We retain for it the same notation, as the domain of $\varphi$ is always clear from the context.
The map $\varphi$  features the following intertwining properties.
\begin{propn}
The isomorphism $\varphi$ intertwines the adjoint actions of $\Ha$ and $\tilde \Ha$
on $\E$.
\end{propn}
\begin{proof}
First of all, remark that  $\Ha$ and $\tilde \Ha$ coincide as associative algebras and enjoy
the same homomorphism $\rho$ to $\E$.
We need to show
\be
\rho(h^{\widetilde{(1)}})\bigl(\varphi (a)\bigr)\rho\bigl(\tilde \gm(h^{\widetilde{(2)}})\bigr)=
\varphi \Bigl(\rho\bigl(h^{(1)}\bigr)a\rho\bigl(\gm(h^{(2)}\bigr)\Bigr)
\label{ad_intertw}
\ee
for all $h\in \Ha$ and $a\in \E$. The twiddled indices on the left-hand side designate the twisted coproduct.
To simplify the formulas we assume $\E=\Ha$ and $\rho=\id$.
For $h\in \Ha$ and $a\in \E$ we find
\be
h^{\widetilde{(1)}}\bigl(\varphi (a)\bigr)\tilde \gm(h^{\widetilde{(2)}})&=&
\F_{1}^{-1}h^{(1)}\F_{1'}\bigl(\F^{-1}_{1''}a\gm(\F^{-1}_{2''})\vt\bigr)\vt^{-1}\gm(\F_{2}^{-1}h^{(2)}\F_{2'})\vt
\nn\\
&=&
\F_{1}^{-1}h^{(1)}a\gm(h^{(2)})\gm(\F_{2}^{-1})\vt=\varphi \Bigl(\rho\bigl(h^{(1)}\bigr)a\rho\bigl(\gm(h^{(2)}\bigr)\Bigr),
\nn
\ee
as required.
\end{proof}
It follows from Theorem \ref{twistadjoint} that
$$\tilde \rho(h):=\rho\bigl(\F_1h\gm(\F_2\zt)\bigr)=\rho\bigl((\F^{-1}_1.h) \F^{-1}_2\bigr)
=(\varphi^{-1}\circ\rho)(h)$$
defines  an algebra  homomorphism $\Ha\to \E_\wr$.
Put $\widetilde{\ad}=\ad_{\tilde \rho}$ to be the adjoint action of $\tilde \Ha$ on $\E_\wr$
defined through the homomorphism $\tilde \rho$:
$$
\widetilde{\ad}(h)a=
\tilde\rho(h^{\widetilde{(1)}})a\tilde\rho\bigl(\tilde \gm(h^{\widetilde{(2)}})\bigr),
 \quad a\in \E_\wr,\quad h\in \tilde \Ha,
\quad h^{\widetilde{(1)}}\tp h^{\widetilde{(2)}}:=\tilde \Delta(h).
$$
The product on the right-hand side is taken with respect to the
multiplication in $\E_\wr$.
Thus one has two representations, $\widetilde{\ad}$ and $\ad$, of the same associative algebra $\Ha$ on the vector space $\E$.
\begin{corollary}
\label{ad_rep}
The representations  $\widetilde{\ad}$ and $\ad$  coincide.
\end{corollary}
\begin{proof}
Apply the homomorphism $\varphi^{-1}\colon \E\to \E_\wr$ to both sides of  equality (\ref{ad_intertw}).
\end{proof}
Theorem \ref{twistadjoint} admits the following generalization.
Consider the smash product $\A\rtimes \Ha$ of the Hopf algebra $\Ha$
and its module algebra $\A$.
The multiplication in $\A\rtimes \Ha$ is described as follows. Both $\A$ and $\Ha$ are subalgebras in
$\A\rtimes \Ha$, which is the free right $\Ha$-module generated by $\A$. Then the multiplication is
determined by the relations $h a=(h^{(1)}\tr a) h^{(2)}$, $h\in \Ha$, $a\in \A$. The subalgebra
$\A$ is invariant under the adjoint action of $\Ha$ on $\A\rtimes \Ha$ due to $\ad(h)(a)=h\tr a$,
for all $h\in \Ha$ and $a\in \A$.
Now consider a twist  $\F\in \Ha\tp \Ha$ and let $\tilde \Ha$ be the twisted Hopf algebra.
Let $\A_\wr$ be the corresponding cotwist of the $\Ha$-module algebra $\A$;
it is a module algebra over
$\tilde\Ha$ with the same  action $\tr$.
\begin{propn}
\label{smash_twisted}
The mapping $\A_\wr\rtimes \tilde \Ha\to \A\rtimes \Ha$,
$a h\mapsto (\F_1\tr a)\F_2h$
is an isomorphism of algebras.
\end{propn}
\begin{proof}
In particular, we need to check that the assignment
$a\mapsto (\F_1\tr a)\F_2$ defines an algebra homomorphism $\A_\wr \to \A\rtimes \Ha$.
The proof is essentially the same as of Theorem \ref{twistadjoint},
because $(\F_1\tr a)\F_2=\ad(\F_1)(a) \>\F_2$. We have yet to check that the cross-relations
in the smash product are preserved.
The equality $h a=(h^{\widetilde{(1)}}\tr a) h^{\widetilde{(2)}}$ in
$\A_\wr\rtimes \tilde \Ha$ goes over to the equality
$
 h(\F_1\tr a) \F_2=(\F_1h^{\widetilde{(1)}}\tr a) \F_2h^{\widetilde{(2)}}
$
in $\A\rtimes \Ha$. The latter equality holds true,
as the right-hand side is nothing but $(h^{(1)}\F_1\tr a)h^{(2)} \F_2$ (follows directly from
the definition of twisted coproduct).
\end{proof}
Remark that Proposition \ref{smash_twisted} may be also considered as a specialization of Theorem \ref{twistadjoint}.
Indeed, the embedding  $\Ha\to\A\rtimes \Ha$
allows to consider $\A\rtimes \Ha$ as an adjoint $\Ha$-module algebra. The factors $\A$ and $\Ha$
are invariant, and the adjoint action restricted to $\A$ is simply $\tr$.
This implies that the  cotwist $(\A\rtimes \Ha)_\wr$ factorizes to the product $\A_\wr\Ha_\wr$, where $\Ha_\wr$ is isomorphic to $\Ha$ as an associative algebra.
This factorization turns into the smash product $\simeq \A_\wr\rtimes \tilde\Ha$ upon
transition from $\Ha_\wr$ to $\Ha$ (it suffices to check the cross-relations only). On the other hand,
$(\A\rtimes \Ha)_\wr$ is isomorphic to $\A\rtimes \Ha$, by Theorem \ref{twistadjoint}.

Adjoint module algebras naturally appear in the following situation.
Suppose $\A$ is an $\Ha$-module algebra and consider
the algebra $\E=\End(\A)$ of linear  endomorphisms of $\A$. The action of $\Ha$ on
$\A$ implies a homomorphism
$\Ha\to \E$. If we define the  adjoint action of $\Ha$ on $\E$
by setting $x.A:=x^{(1)}A\gm(x^{(2)})$, then the action of $\E$ on $\A$ is $\Ha$-equivariant
(that is, the map $\E\tp \A \to \A$ is a homomorphism of $\Ha$-modules).
Twist of $\Ha$ induces a transformation of the algebras $\A$ and  $\E$, as well as of the $\E$-action on $\A$:
$$
X*a:=(\F_1.X)(\F_2 a),
\quad X\in \E, \quad a\in \A.
$$
The algebra isomorphism $\varphi^{-1}\colon \E\to \E_\wr$ and the transformation of the action compensate each other in the following sense:
$$
\varphi^{-1}(X)*a=X a
$$
for all $X\in \E$ and $a\in \A$.

\section{$\Ha$-Lie algebras}

In the present section, we recall braided or $\Ha$-Lie algebras in our terminology, assuming $\Ha$ to be triangular.

Suppose that $\Ha$ is triangular and let $\g$ be an $\Ha$-module. Let us call $\g$ an
$\Ha$-Lie algebra if it is equipped with an $\Ha$-equivariant map $[.,.]_\Ru\colon\g\tp \g\to \g$
satisfying the braided Jacobi identity, \cite{G},
\be
\label{Jacobi}
[\xi,[\eta,\zt]_\Ru]_\Ru=[[\xi,\eta]_\Ru,\zt]_\Ru+[\Ru_2\tr\eta,[\Ru_1\tr\xi,\zt]_\Ru]_\Ru.
\ee
plus the skew symmetry condition
$$
[\Ru_2\tr \xi,\Ru_1\tr \eta]_\Ru=-[\eta,\xi]_\Ru.
$$
For example,
any associative $\Ha$-module algebra $\B$
is an $\Ha$-Lie algebra via the commutator
$$
[a,b]_\Ru:=ab-(\Ru_2\tr b) (\Ru_1\tr a)
$$
for all $a,b\in \B$. The algebra $\B$ is quasi-commutative if and only if this commutator vanishes.
In general, we shall call an $\Ha$-Lie algebra with zero commutator quasi-Abelian.

One can naturally define homomorphisms between $\Ha$-Lie algebras, which are equivariant maps preserving
the commutators.

The quotient of the tensor algebra $\T(\g)$ by the relations
$\xi\tp \eta-(\Ru_2\tr \eta)\tp  (\Ru_1\tr \xi)-[\xi,\eta]_\Ru$
is called universal enveloping algebra of $\g$
and denoted by $\Ur(\g)$. Since the ideal is invariant, the $\Ha$-action on
$\T(\g)$ produces an action on $\Ur(\g)$.

\begin{propn}
Suppose that $\Ha$ is triangular and let $\g$ be an $\Ha$-Lie algebra.
Then the smash product $\Ur(\g)\rtimes \Ha$ is a triangular Hopf algebra with the universal R-matrix $\Ru$ and the
comultiplication extended from $\Ha$ by
\be
\label{coprod}
\Delta(\xi):=\xi\tp 1+\Ru_2\tp \Ru_1\tr\xi, \quad \xi\in \g.
\ee
The counit and antipode when restricted to $\g$  are given by
$\ve(\xi)=0$ and $\gm(\xi)=-(\Ru_2\tr\xi)\Ru_1$.
\label{uni_H_smash}
\end{propn}
\begin{proof}
The formulas for counit and antipode readily follow from the comultiplication, which is obviously coassociative.
Let us prove that it is a homomorphism.
First of all, the comultiplication respects the cross-relations $h\xi=(h^{(1)}\tr\xi) h^{(2)}$:
\be
\Delta(h)\Delta(\xi)&=&h^{(1)}\xi\tp h^{(2)}+h^{(1)}\Ru_2\tp h^{(2)}\Ru_1\tr\xi
\nn
\\&=&
(h^{(1)}\tr\xi)h^{(2)}\tp h^{(3)}+h^{(1)}\Ru_2\tp (h^{(2)}\Ru_1)\tr\xi h^{(3)}
\label{cross_lie_H}\\&=&
(h^{(1)}\tr\xi)h^{(2)}\tp h^{(3)}+\Ru_2h^{(2)}\tp \Ru_1\tr(h^{(1)}\tr\xi) h^{(3)}
=\Delta(h^{(1)}\tr\xi)\Delta(h^{(2)}),
\nn
\ee
for all $\xi\in \g$, $h\in \Ha$.
Further, for any tensor $\xi\tp \eta\in \g\tp\g$ we write
\be
\Delta(\xi\eta)&=&
\xi\eta\tp 1+\Ru_2\tp \Ru_1\tr(\xi\eta)
+\xi\Ru_2\tp \Ru_1\tr\eta
+\Ru_2\eta\tp \Ru_1\tr\xi,\nn\\
\Delta(\tilde\eta\tilde\xi)&=&
\tilde\eta\tilde\xi\tp 1+\Ru_2\tp \Ru_1\tr(\tilde\eta\tilde\xi)
+\tilde\eta\Ru_2\tp \Ru_1\tr\tilde\xi
+\Ru_2\tilde\xi\tp \Ru_1\tr\tilde\eta .\nn
\ee
where we put $\tilde \eta\tp \tilde \xi= \Ru_2\tr\eta\tp \Ru_1\tr\xi=
\Phi_1\eta \Phi_2\tp \Phi_3\xi \Phi_4$ for certain  $\Phi=\Phi_1\tp \Phi_2 \tp \Phi_3\tp \Phi_4\in \Ha^{\tp 4}$.
Subtracting $ \Delta(\tilde\eta\tilde\xi)$ from $\Delta(\xi\eta)$ we prove the statement,  provided that
$$
\Ru_{2}\eta \tp (\Ru_{1}\tr\xi)
=
(\Phi_1\eta\Phi_2)\Ru_2 \tp  \Ru_1 \tr(\Phi_3\xi\Phi_4),
$$
$$
\xi\Ru_2 \tp  (\Ru_1 \tr \eta)
=
\Ru_{2}(\Phi_3 \xi \Phi_4)\tp  \Ru_{1}\tr(\Phi_1\eta \Phi_2).
$$
They are satisfied indeed, because, by the definition of $\Phi$,
$$
(\Delta\tp \Delta^{op})(\Ru)
=
\gm^{-1}_2\gm_3\Phi.
$$
Then the first equation holds identically.
The second equation is true for triangular $\Ru$.
\end{proof}
There is a natural $\Ha$-invariant filtration on $\Ur(\g)$, and the corresponding graded algebra is
the quasi-commutative $\Ru$-symmetric algebra of the module $\g$. Using this filtration, one can check that $\g$ is
precisely the set of quasi-primitive elements in $\Ur(\g)$, i.e. the elements satisfying
the identity (\ref{coprod}).

There is a natural action of the Hopf algebra  $\Ur(\g)\rtimes \Ha$ on  $\Ur(\g)$ making the latter
an  $\Ur(\g)\rtimes \Ha$-module algebra. Namely, $\Ur(\g)$ is invariant under the adjoint action.
So, for any $\xi\in \g$  and $a\in \Ur(\g)$ one has $\ad(\xi)a=[\xi,a]_\Ru$.

\begin{propn}
Let $\F\in \Ha\tp \Ha$ be a twisting cocycle.
Suppose that $\g$ is an $\Ha$-Lie algebra.
Then the
commutator $[\>.\>,\>.\>]_{\tilde \Ru}:=[\>.\>,\>.\>]_{\Ru}\circ \F$
defines on $\g$ a structure of an $\tilde\Ha$-Lie algebra.
\end{propn}
\begin{proof}
We need to show that
$[\>.\>,\>.\>]_{\tilde \Ru}$ satisfies the Jacobi identity.
Denote by $\tau$ the ordinary flip $\g\tp\g\to \g\tp \g$ and
write for the mapping $[\>.\>,[\>.\>,\>.\>]_{\tilde \Ru}]_{\tilde \Ru}\colon\g\tp\g\tp\g\to \g$:
$$
[\>.\>,[\>.\>,\>.\>]_{ \Ru}]_{\Ru}\circ \Delta_2(\F)\F_{23}
=[[\>.\>,\>.\>]_{ \Ru},\>.\>]_{\Ru}\circ \Delta_1(\F)\F_{12}+
[\>.\>,[\>.\>,\>.\>]_{ \Ru}]_{\Ru}\circ \tau_{12}\circ\Ru_{12} \Delta_1(\F)\F_{12}.
$$
Here we applied the twist equation and the Jacobi identity to the commutator $[\>.\>,\>.\>]_{\Ru}$.
The first term on the right-hand side is just $[[\>.\>,\>.\>]_{ \tilde\Ru},\>.\>]_{\tilde\Ru}$.
Using the standard Hopf algebra technique, we get
$$
\tau_{12}\circ\Ru_{12} \Delta_1(\F)\F_{12}
=\tau_{12}\circ\Delta^{op}_1(\F)\F_{21}\tilde \Ru_{12}
=\Delta_1(\F)\F_{12}\circ\tau_{12}\circ\tilde \Ru_{12},
$$
which gives for the second term the expression $[\>.\>,[\>.\>,\>.\>]_{\tilde \Ru}]_{\tilde \Ru}\circ\tau_{12}\circ\tilde \Ru_{12}$, i.e. the
second summand in the Jacobi identity.
\end{proof}

We denote by $\g_\wr$ the vector space $\g$ together with
this $\tilde \Ha$-Lie algebra structure.
\begin{propn}
The cotwist  ${\Ur(\g)}_\wr$ is isomorphic to  $\Ur(\g_\wr)$,
and the isomorphism is extended from the identity map $\g \ni\xi\mapsto \xi\in\g_\wr$.
\label{Utw}
\end{propn}
\begin{proof}
The universal enveloping algebra $\Ur(\g)$ is a quotient of the tensor algebra $\T(\g)$
over the ideal $J(\g)$ generated by the commutator relations. It is easy to see that the cotwist $\Ur(\g)_\wr$
is isomorphic to  the quotient of  $\T(\g)$ by the ideal $J(\g_\wr)$,
where $\g_\wr$ is obtained from $\g$ by rewriting the defining relations in
terms of the new multiplication, $*$.
That is, $(\F^{-1}_1\tr\xi)*(\F^{-1}_2\tr\eta)-(\F^{-1}_1\Ru_2\tr\eta)*(\F^{-1}_2\Ru_1\tr\xi)=[\xi,\eta]_{\Ru}$
for all $\xi,\eta\in \g$. As $\F$ is invertible, this is equivalent
to the relations
$\xi*\eta-(\tilde\Ru_2\tr\eta)*(\tilde\Ru_1\tr\xi)=[\F_1\tr\xi,\F_2\tr\eta]_{\Ru}$,
as required.

The isomorphism $\Ur(\g_\wr)\to \Ur(\g)_\wr$
is induced by a linear endomorphism
of the tensor algebra $T(\g)$. On the homogeneous component of degree $k$ this isomorphism
is given by $\xi_1\tp \ldots \tp \xi_k\mapsto \F^{(k)}_1\tr\xi_1\tp  \ldots \tp \F^{(k)}_k\tr\xi_k$,
where $\F^{(k)}\in \Ha^{\tp k}$ is defined inductively
by $\F^{(1)}:=1$, $\F^{(2)}:=\F$,  $\F^{(m+n)}:=(\Delta^{(m)}\tp\Delta^{(n)})(\F)(\F^{(m)}\tp \F^{(n)})$.
This definition is consistent, i.e., independent of the partition $k=m+n$.
Thus the isomorphism in question is identical on $\g$ by construction.
\end{proof}

\begin{propn}
\label{smash_twisted_1}
Consider a twisting cocycle $\F\in \Ha\tp \Ha$ as that of the Hopf algebra
$\Ur(\g)\rtimes \Ha$.  Then the assignment
$\g\ni \xi\mapsto (\F_1\tr \xi)\F_2\in \Ur(\g)\rtimes \Ha$
defines an isomorphism of Hopf algebras
\be
\Ur(\g_\wr)\rtimes \tilde \Ha\to \widetilde{\Ur(\g)\rtimes \Ha}
\label{smash_uni}
\ee
which is identical on $\Ha$.
\end{propn}
\begin{proof}
It follows from Proposition \ref{smash_twisted} that
the assignment $\varphi\colon \xi\mapsto(\F_1\tr \xi)\F_2$ for all $\xi\in \g$ is a restriction of the associative algebra isomorphism
$\Ur(\g)_\wr\rtimes \tilde \Ha\to \Ur(\g)\rtimes \Ha$.
On the other hand, the identity mapping $\g_\wr\to \g$ extends to
an $\tilde \Ha$-equivariant isomorphism $\Ur(\g_\wr)\to \Ur(\g)_\wr$. Therefore
the assignment in question extends to an isomorphism of algebras.

We need to check that the coalgebra structure is preserved by (\ref{smash_uni}). That is obvious for its
restriction to $\tilde \Ha$. Let as prove that $\varphi$ is a coalgebra map when restricted
to $\Ur(\g_\wr)$. To this end, we calculate
$
\tilde \Delta\bigl(\varphi(\xi)\bigr)=\F^{-1}\Delta\bigl((\F_1\tr\xi) \F_2\bigr)\F
$
for $\xi\in \g_\wr$. Note that all the products and coproducts here are written in terms of
$\Ur(\g)\rtimes \Ha$.
Using  the twist identity we get for
$
\tilde \Delta\bigl(\varphi(\xi)\bigr)
$:
\be
&&\F^{-1}\Delta(\F_1\tr\xi) (\F_2^{(1)}\F_{1'}\tp\F_2^{(2)}\F_{2'})
=\F^{-1}\Delta\bigl((\F_1^{(1)}\F_{1'})\tr\xi \bigr)(\F_1^{(2)}\F_{2'}\tp\F_2 )\nn.
\ee
Now evaluate $\Delta$ of the element from $\g$ and get the sum of two terms, of which the first is
the product of $\F^{-1}$ and $\bigl((\F_1^{(1)}\F_{1'})\tr\xi\bigr) \F_1^{(2)}\F_{2'}\tp\F_2 =\F_1(\F_{1'})\tr\xi\F_{2'}\tp\F_2$. Thus, the first summand is equal to $\varphi(\xi)\tp 1$.
The other summand is the product of $\F^{-1}$ and
\be
&&
\Ru_2\F_1^{(2)}\F_{2'}\tp\bigl((\Ru_1\F_1^{(1)}\F_{1'})\tr\xi\bigr) \F_2 =
\F_1^{(1)}\Ru_2\F_{2'}\tp\bigl((\F_1^{(2)}\Ru_1\F_{1'})\tr\xi\bigr) \F_2
\nn\\
&&=
\F_1^{(1)}\F_{1'}\tilde\Ru_2 \tp\bigl((\F_1^{(2)}\F_{2'}\tilde \Ru_1)\tr\xi\bigr) \F_2
=
\F_1\tilde\Ru_2 \tp\bigl((\F_2^{(1)}\F_{1'}\tilde \Ru_1)\tr\xi \F_2^{(2)}\bigr)\F_{2'}
\nn\\
&&=
\F_1\tilde\Ru_2 \tp\F_2\bigl((\F_{1'}\tilde \Ru_1)\tr\xi \bigr)\F_{2'}.
 \nn
\ee
Multiplying this by $\F^{-1}$ we get $\tilde\Ru_2 \tp\varphi(\tilde \Ru_1\tr\xi)$ for the second
summand. Note that the dot here is the adjoint action of $\Ha$ on $\g$, not of $\tilde \Ha$.
Apply to this formula (\ref{ad_intertw}) where $\rho$ is the embedding $\Ha\hookrightarrow \Ur(\g)\rtimes \Ha$ and get $\tilde\Ru_2 \tp\tilde \Ru_1\tr\varphi(\xi)$. This proves that $\varphi$ respects comultiplication
and when restricted to $\g_\wr$ is therefore a coalgebra map.
\end{proof}
\begin{remark}
\label{g-g_cotw}
Let us emphasize that the $\tilde \Ha$-Lie algebra $\g_\wr $ is rotated from $\g$ inside of the associative algebra $\Ur(\g)\rtimes \Ha$
by means of $\varphi$, that is $\g_\wr =\varphi(\g)$ .
This  readily follows from the proof of the above proposition.
\end{remark}

For any triangular Hopf algebra $\Ha$ denote by $\Lie(\Ha)$
the subset of elements satisfying the condition
$$
\Delta(\xi):=\xi\tp 1+\Ru_2\tp \Ru_1\tr\xi=\xi\tp 1+\Ru_2\Ru_{1'}\tp \Ru_1\xi\Ru_{2'}.
$$
It is easy to check that $\Lie(\Ha)$ is invariant under the adjoint action of $\Ha$.
Moreover,  $\Lie(\Ha)$ is an $\Ha$-Lie algebra under the $\Ru$-commutator,
hence there exists an algebra homomorphism $\Ur\bigl(\Lie(\Ha)\bigr)\to \Ha$.
This homomorphism also extends to a Hopf algebra homomorphism $\Ur\bigl(\Lie(\Ha)\bigr)\rtimes \Ha\to \Ha$.

\begin{propn}
\label{LHtw}
Put $\g=\Lie(\Ha)$ and denote by $\g_\wr$ the cotwist of $\g$ by a cocycle $\F\in \Ha\tp\Ha$.
Then $\Lie(\tilde\Ha)\simeq \g_\wr$.
\end{propn}
\begin{proof}
Denote by $\tilde \g$ the Lie algebra $\Lie (\tilde\Ha)$.
Apply the chain of Hopf algebra homomorphisms
$$
\Ur(\g_\wr)\rtimes \tilde \Ha\to \widetilde{\Ur(\g)\rtimes \Ha} \to \tilde \Ha
$$
to $\g_\wr$ and find that $\varphi(\g_\wr)\subset \tilde \g$, see Remark \ref{g-g_cotw}.
On the other hand, twist is invertible
and the role of $\varphi$ for the inverse twist belongs to $\varphi^{-1}$. Hence we can do the same with $\tilde \Ha$ replaced by $\Ha$ and get the inverse inclusion
$\varphi^{-1}(\tilde \g_\wr)\subset \g$.  This implies the statement, because $\g_\wr$ is isomorphic to $\g$ as a vector space.
\end{proof}

\begin{propn}
Suppose that $\Ha$ is generated by $\Lie(\Ha)$ as an associative algebra.
Then $\tilde \Ha$ is generated by $\Lie(\tilde\Ha)$.
\end{propn}
\begin{proof}
Denote $\g=\Lie(\Ha)$. It is sufficient to consider the case when $\Ha=\Ur(\g)$.
Put $\A:=\Ha$ and consider it as the adjoint module algebra over $\Ha$.
Since $\A=\Ur(\g)$, one has $\tilde \A=\Ur(\g)_\wr\simeq\Ur(\g_\wr)$,
by Proposition \ref{Utw} and  $\Lie(\tilde\Ha)\simeq\g_\wr$  by Proposition \ref{LHtw}.
But $\tilde\A$ is isomorphic to $\A$, by Proposition \ref{twistadjoint}.
\end{proof}
Note that the Lie algebras $\g_\wr$ and $\g$ are related as subsets of $\Ha$ by the isomorphism $\varphi$:
$\g_\wr=\varphi(\g)$, cf. Remark \ref{g-g_cotw}.

Suppose that $\Ha$ is a classical universal enveloping algebra $\Ur(\g)$.
Given a twist  of $\Ha$, the associative algebra $\tilde \Ha$ is isomorphic to  $\Ur(\g_\wr)$.
The elements of $\g_\wr$ form a Lie algebra over $\tilde \Ha$, and $\g_\wr$
is linked to $\g$ by the transformation $\varphi$.
In other words, there is a set of generators for $\Ur(\g)$, such
that the twisted comultiplication is given
by $\tilde \Delta(\xi)=\xi\tp 1+\Ru^{-1}_1\Ru_{1'}\tp \Ru^{-1}_2 \xi \Ru_{2'}$,
where $\Ru=\F_{21}^{-1}\F$.

\section{Twisted vector fields}
In this section we consider an important example of quantum Lie algebras formed by "quantum vector fields", \cite{Wor}.
As before, we assume  that $\Ha$ is a triangular Hopf algebra and $\A$ is a left $\Ha$-module algebra.

\begin{definition}
An element $X\in \End(\A)$  is called a first order differential operator
if
$$
X(ab)=(Xa)b+(\Ru_{2}\Ru^{-1}_{2'}a) (\Ru_1 X\Ru^{-1}_{1'}b)=(Xa)b+(\Ru_{2}a) (\Ru_1. X b).
$$
\end{definition}
The subset of first order differential operators is denoted by $\Der_\Ha(\A)$. Following the geometric analogy,
we shall also call elements of $\Der_\Ha(\A)$ vector fields, thinking of  $\A$ as the function algebra
on a quantum space. If $\A$ is quasi-commutative, then $\Der_\Ha(\A)$ is a natural $\A$-submodule
in the left module $\End(\A)$. Clearly the $\A$-action on $\Der_\Ha(\A)$ is $\Ha$-equivariant.

We can introduce a "comultiplication" on vector fields as a map
$$
\Delta\colon\Der_\Ha(\A)\to  \Der_\Ha(\A)\rtimes \Ha\tp \Der_\Ha(\A),
\quad
\Delta \colon X\mapsto X\tp 1 + \Ru_{2} \tp \Ru_1. X .
$$
This allows to use same technique to prove facts about $\Der_\Ha(\A)$ as we did in the previous section for
$\Ha$-Lie algebras.
\begin{propn}
The set $\Der_\Ha(\A)$ is $\Ha$-invariant.
\end{propn}
\begin{proof}
Same as the proof of (\ref{cross_lie_H}).
\end{proof}
If the algebra $\A$ is quasi-commutative, $\Der_\Ha(\A)$ is a left $\A$-module regarded as a subalgebra of $\End(\A)$ by left multiplication.

Now consider the cotwist of $\A$, and also of  $\E:=\End(\A)$.
The new action of $\E_\wr$ on $\A_\wr$ is expressed through the old action
and the homomorphism $\varphi$ by
$X*a=\F_1.X \F_2 a=\varphi(X)a$.


\begin{propn}
\label{DerisFinv}
$\Der_{\tilde \Ha}(\A_\wr)=\Der_\Ha(\A)$ as subsets in $\E$.
\end{propn}
\begin{proof}
Suppose $X\in \Der_\Ha(\A)$. We need to show
that
$$
X*(a*b)=\bigl(X*a\bigr)* b+
(\tilde\Ru_2 a)*\bigl((\tilde \Ru_1.X )* b\bigr)
$$
for all $a,b\in \A$.
Applying to the left-hand side the same steps as in the proof of Proposition \ref{smash_twisted_1} we arrive
at the expression
$
\bigl(X*a\bigr)* b+
(\tilde\Ru_2 a)*\bigl((\F_{1'}\tilde \Ru_1).X \F_{2'} b\bigr).
$
The dot in the second term involves the coproduct $\Ha$.
The intertwining formula (\ref{ad_intertw}) allows to write it as
$(\tilde\Ru_2\tilde\Ru_{1'} a)*\bigl(\tilde \Ru_1. X * \tilde\Ru_{2'} b\bigr)$, i.e. to interpret the dot
trough the new coproduct.
\end{proof}
The following statement asserts that
the subset $\Der_\Ha(\A)\subset \End(\A)$ is a subalgebra of the commutator
$\Ha$-Lie algebra $\End(\A)$.

\begin{propn}
\label{Derisclosed}
For any  $X,Y\in  \Der_\Ha(\A)$,
the operator
$
[X,Y]_{\Ru}=XY-(\Ru_2.Y)(\Ru_1.X)
$
 belongs to $\Der_\Ha(\A)$.
\end{propn}
\begin{proof}
Same as the proof of Proposition \ref{uni_H_smash}.
\end{proof}

Thus, the $\Ha$-module $\g=\Der_\Ha(\A)$ is an $\Ha$-Lie algebra. The action of $\g$ on $\A$ extends to an
action of the universal enveloping algebra $\Ur(\g)$, which together with the action of $\Ha$ makes $\A$
a module algebra over $\Ur(\g)\rtimes \Ha$.
The representation $\Ha\to \End(\A)$ induces an $\Ha$-Lie algebra homomorphism
$\Lie(\Ha)\to \Der_\Ha(\A)$.

Further we give examples of $\Ha$-Lie algebras naturally arising in geometrical applications
and, in particular, in gauge field theory on non-commutative space-time.
\begin{example}
Let $\Ha$ be a triangular Hopf algebra and $\A$ be a quasi-commutative $\Ha$-module algebra.
Suppose that $\g$ is an $\Ha$-Lie algebra.
Denote by $\A\bowtie \g$ the $\Ha$-module $\A\tp \g$ equipped with the bilinear operation
\be
[a\tp \xi,b\tp \eta]_{\Ru}:=a(\Ru_2\tr b)\tp [\Ru_1\tr \xi,\eta]_{\Ru},
\label{current_commutator}
\ee
where $a,b\in \A$ and $\xi,\eta\in \g$.
\begin{propn}
\label{current}
The $\Ha$-module $\A\bowtie \g$ is an $\Ha$-Lie algebra.
\end{propn}
\begin{proof}
Denote by the $\A\bowtie \Ur(\g)$ the associative algebra with the underlining vector space $\A\tp \g$
and the multiplication
$
(a\tp \xi)(b\tp \eta):=a(\Ru_2\tr b)\tp (\Ru_1\tr \xi)\eta.
$
This is the so called braided tensor product of $\A$ and $\Ur(\g)$ and is  a module algebra over $\Ha$.
The vector space $\A\bowtie \g$ is a $\Ha$ -submodule in $\A\bowtie \Ur(\g)$.
Let us prove that the commutator $\Ha$-Lie algebra restricts to $\A\bowtie \g$ and coincides with
the pre-defined $\Ha$-Lie algebra structure on $\A\bowtie \g$. That is obvious for
the $\Ha$-Lie subalgebra $\g\subset \A\rtimes \g$ and follows from the very construction. Further, every element $\xi\in \g$
$\Ru$-commutes with every element $b\in \A\subset \A\bowtie \Ur(\g)$, hence
the commutator of $\xi$ and $b\tp \eta$ has the desired form. By assumption, the algebra $\A$ is quasi-commutative.
From this we conclude that derivations of $\A\bowtie \Ur(\g)$ form a left $\A$-module
(in this particular case even two-sided module). Hence the commutator of $a\tp \xi$ and $b\tp \eta$
is given by (\ref{current_commutator}).
\end{proof}
The algebra $\A\bowtie \g$ is an $\A$-bimodule and
the commutator is a bimodule map.
The left action comes from the regular $\A$-action on itself, while
the right action is defined by $(a\tp \xi)b=a\Ru_2\tr b\tp \Ru_1\tr\xi$.
Geometrically, such Lie algebras  are modeled by vector bundles whose fiber is an $\Ha$-Lie algebra,
say, of the gauge group.

\end{example}
\begin{example}
The previous example is a special case of the following construction.
Suppose that an $\Ha$-Lie algebra $\g$ acts on  a quasi-commutative $\Ha$-module algebra $\A$.
That means that there is a ($\Ha$-equivariant) homomorphism $\g\to \Der_\Ha(\A)$.
Regard the tensor product $\A\tp \g$ as the natural left $\A$-module and identify $\g$ with
$1\tp \g \subset \A\tp \g$.
There is a unique extension of the $\Ha$-Lie algebra structure from $\g$ to $\A\tp \g$
satisfying
$$
[\xi, f\eta]_{\Ru}=(\xi f)\eta+\Ru_2f [\Ru_1\xi, \eta]_{\Ru}, \quad \forall \xi,\eta \in \A\tp \g, \forall f\in \A.
$$
This Lie algebra may be regarded as that of local transformations, as opposed to the global transformations
by $\g$.
The loop algebra from the previous example is obtained from this by taking the zero action of $\g$ on $\A$.
\end{example}
\section*{Acknowledgements}
This research is supported in part by the RFBR grant 09-01-00504.

\end{document}